\input amstex
\documentstyle{amsppt}
\topmatter
\title Projective Resolutions for Graph Products
\endtitle
\author Daniel E. Cohen
\endauthor
\affil Queen Mary and Westfield College, London University
\endaffil
\address School of Mathematical Sciences, Queen Mary and Westfield 
College, Mile End Road, London E1 4NS, England
\endaddress
\email D.E.Cohen\@uk.ac.qmw.maths
\endemail

\abstract 

Let $\Gamma$ be a finite graph together with a group $G_v$ at each vertex 
$v$. The {\it graph product\/} $G(\Gamma)$ is obtained from the free product 
of all $G_v$ by factoring out by the normal subgroup generated by $\{ 
g^{-1}h^{-1}gh; g \in G_v, h \in G_w \}$ for all adjacent $v,w$.

In this note we construct a projective resolution for $G(\Gamma)$ given 
projective resolutions for each $G_v$, and obtain some applications.
\endabstract 
\endtopmatter
\document
Let $\Gamma$ be a finite graph together with a group $G_v$ at each vertex 
$v$. The {\it graph product\/} $G(\Gamma)$ is obtained from the free product 
of all $G_v$ by factoring out by the normal subgroup generated by $\{ 
g^{-1}h^{-1}gh; g \in G_v, h \in G_w \}$ for all adjacent $v,w$.

In this note we construct a projective resolution for $G(\Gamma)$ given 
projective resolutions for each $G_v$, and obtain some applications. This is 
quite easy to do, since $G(\Gamma)$ is built up from the vertex groups by 
direct products and amalgamated free products.

Let $A$ be an arbitrary group, and let $R$ be a commutative ring. A {\it 
projective resolution\/} for $A$ is an exact sequence $$ \text{\bf P}: \ldots 
\to P_n \overset d_n\to\rightarrow P_{n-1} \to \ldots P_1 \to P_0 \to R \to 
0$$ of projective (right) $R A$\/-modules. We shall always assume that $P_0 = 
R A$.

Let $G = A *_C B$ and let {\bf Q} and {\bf N} be projective resolutions 
of $B$ and $C$. Suppose that $\text{\bf N}\otimes_{R C } R A$ and $\text{\bf 
N} \otimes_{R C }R B$ are summands of {\bf P} and {\bf Q} respectively. Let 
$S_k$ be the quotient of $(P_k \otimes_{R A }R G) \oplus (Q_k \otimes_{R  B} R
G)  $ by the submodule generated by all $(n\otimes_{R C} 1, -n\otimes_{R C} 1)$, 
with the obvious map from $S_k$ to $S_{k-1}$.

\proclaim{Lemma 1} In the above situation, $$\ldots S_k \to S_{k-1} \to 
\ldots \to S_1 \to R G \to R \to 0$$ is a projective resolution for $G$. 
\endproclaim

{\bf Proof}. It is easy to check that each $S_k$ is projective.

Let $I_A$ be the augmentation ideal of $A$, with similar notation for 
$B,C,G$. We can regard the resolution {\bf P} as finishing with $P_2 \to P_1 
\to I_A \to 0$, and this sequence is exact. It is then immediate, by 
elementary properties of exact sequences, that we have an exact sequence $$ 
\ldots S_n \to S_{n-1} \to \ldots \to S_2 \to S_1 \to M \to 0,$$ where $M$ 
is defined to be the quotient of $(I_A \otimes_{R  A}R G) \oplus (I_B \otimes_{R 
B} R G)$ by the submodule generated by all $((c-1)\otimes_{R C}1, - 
(c-1)\otimes_{R C}1)$.

However, according to Lemma 4.3 and Theorem 4.7 of [6], $M$ is exactly 
$I_G$, as required. \qed

Now let $H=A \times B$. We find a projective resolution for $H$ as the 
tensor product of {\bf P} and {\bf Q}. Precisely, we define $(\text{{\bf P}} 
{\otimes ^ +} \text{{\bf Q}})_n$ to be $\sum_{i+j=n}P_i \otimes_R  Q_j$, where 
$i$ and $j$ are non-zero. The superscript $+$ is used because the tensor 
product is defined to allow $i$ and $j$ to be zero. Our modification will, in 
fact, lead to a messy formula for the resolution of $H$. It is used because 
at a later point it ensures that we do not have summands repeated 
unnecessarily.

Note that $(\text{{\bf P}} {\otimes ^ +} \text{{\bf Q}})_n$ is an {\it 
H\/}-module, as are $P_n \otimes_{R A} R H$ and $Q_n \otimes_{R B} R H$. 
These are plainly projective modules. Further, these last two modules are 
isomorphic to $P_n \otimes_R R B$ and $R A \otimes_R Q_n$, and we identify 
the isomorphic modules.

Let {\bf T} be $(\text{{\bf P}} {\otimes ^+}  \text{{\bf Q}}) \oplus 
(\text{{\bf P}}\otimes_{R A} R H) \oplus (\text{{\bf Q}}\otimes_{R B}R H)$ in 
dimensions greater than 0, while $T_0=R H$. We define a boundary operator on 
{\bf T} in the obvious way on the last two summands. On $P_i \otimes_R Q_j$ 
we define it to map $p\otimes q$ to $(pd_i)\otimes q + (-1)^i p \otimes 
(qd_j)$. When $i$ or $j$ is 1, one (or both) of these terms lies in the third 
or second of the summands.

As is well-known (for instance, V.1.1 of [3]) we have the following 
lemma.

\proclaim{Lemma 2} {\bf T} is a  projective
resolution for $H$.  \endproclaim

Now let $\Gamma$ be a finite graph with vertex group $G_v$ at the vertex 
$v$. Let $\Delta$ be any full subgraph of $\Gamma$, with the corresponding 
graph product $G(\Delta)$. There is an obvious homomorphism from $G(\Delta)$ 
to $G(\Gamma)$ induced by inclusion of the vertex groups, and another 
homomorphism from $G(\Gamma)$ to $G(\Delta)$ induced by the inclusion on 
those vertex groups with vertices in $\Delta$ and trivial on the other vertex 
groups. It follows that we have a monomorphism, which we regard as inclusion, 
from $G(\Delta)$ to $G(\Gamma)$.

Let $v$ be any vertex, let $\Delta$ be the full subgraph on all vertices 
except $v$, and let $\Omega$ be the full subgraph on $v$ and all vertices 
adjacent to to $v$. Then $G(\Gamma) = (G_v \times G(\Omega \cap \Delta) 
*_{G(\Omega \cap \Delta)} G(\Delta)$. This was proved by Green [7], and it 
is easy to obtain homomorphisms between these two groups which are inverses 
of each other. From this, and our lemmas, it is easy to find a projective 
resolution for $G(\Gamma)$, but we need some more notation. We let $\text{\bf 
P}^ v$ be a projective resolution of $G_v$, with dimension 0 being omitted 
for later notational convenience.

It is convenient to regard the vertices as having a specified total 
order. Let $K$ be a complete subgraph of $\Gamma$ whose vertex set (in the 
given order) is $\{v_1,\ldots ,v_k\}$. Let $\text{{\bf P}}^ K = (\text{{\bf
P}}^  {v_1} \otimes ^ +_R \ldots \otimes ^ + _R \text{{\bf P}}^ {v_k})\otimes_{R 
G(K)} G(\Gamma)$. As in the case of a direct product (which corresponds to 
the complete graph on two vertices), a boundary operator is defined on $\sum 
\text{{\bf P}}^ L$, where the sum is taken over all complete subgraphs $L$ of 
$K$.

We note that there is a bijection between the complete subgraphs of 
$\Omega \cap \Delta$ and the complete subgraphs of $\Omega$ which contain $v$ 
and at least one other vertex. We will choose $v$ to be the first vertex of 
the ordering.

Using the two lemmas, the theorem below now follows easily by induction 
on the number of vertices of $\Gamma$. 

\proclaim{Theorem} $\sum \text{{\bf P}}^ K$, taken over all complete 
subgraphs $K$ of $\Gamma$, forms the positive-dimensional part of a 
projective resolution for $G(\Gamma)$. \endproclaim

Recall that a group is $FP_n$ ($FP_{\infty}$) if it has a projective 
resolution finitely generated in dimensions $\le n$ (in all dimensions), it 
is of finite cohomological dimension if it has a projective resolution zero 
in all but finitely many dimensions, and it is $FP$ if it has a projective 
resolution finitely generated in all dimensions and zero in all but finitely 
many. The following corollary is immediate from the theorem.

\proclaim{Corollary 1} Let the vertex groups of a graph product over a 
finite graph all be $FP_n$ (for any $n \le \infty$), of finite cohomological 
dimension, or $FP$. Then the graph product has the same property. 
\endproclaim

Chiswell [4] defines an Euler characteristic for any $FP$ group, and he 
obtains in [5] a formula for the Euler characteristic of a graph product, 
using the inductive decomposition of a graph product in terms of amalgamated 
free products and direct products. Our theorem (which relies on this 
decompostion) makes it clearer where Chiswell's formula comes from. However, 
his proof also holds for Brown's definition [3] of the Euler characteristic, 
which ours does not appear to.

Let $G$ be a group, and let $P$ be a finitely generated projective 
$RG$-module. Then Hattori [8] and Stallings [9] define a trace function $t_P$,
which  is an element of $R(G/G')$. Let $\tau _P$ be the coefficient of 1 in
$t_P$.  When $G$ is $FP$ take a projective resolution {\bf P} which is finitely 
generated in all dimensions and zero in all but finitely many. Chiswell 
defines the Euler characteristic of $G$ to be $\sum_0 ^ {\infty} \tau 
_{P_i}$.

Let $G$ be a subgroup of $H$, and let $P$ be a finitely generated 
projective $RG$-module. Then $P \otimes _{RG} RH$ is a finitely generated 
projective $RH$-module, and $t_{P \otimes _{RG} RH}$ is the image in $R(H/H') $
of $t_P$; hence $\tau _ {P \otimes _{RG} RH} = \tau _P$. Also, if $Q$ is a 
finitely generated projective $RF$-module, for some group $F$, then $P 
\otimes _R Q$ is a finitely generated projective ($G \times F$)-module, and 
$t_{P \times _R Q} = t_P \otimes t_Q$ in $R(G \times F)/ (G \times F)' = 
R(G/G') \otimes _R R(F/F')$; hence $\tau_{P \otimes _R Q} = \tau _P \tau _Q$. 
These facts can be found in [2], and are easy to prove directly.

Let $\Gamma$ be a finite graph, with an $FP$ group at each vertex. For 
each complete subgraph $K$ of $\Gamma$, with vertices $v_1, \ldots, v_k$, let 
$\chi _K =(\chi _{v_1} -1)\ldots (\chi_ {v_k} -1)$, where $\chi_{v_i}$ is the 
Euler characteristic of the group at the vertex $v_i$. From our main theorem 
and the remarks in the previous paragraph, we immediately deduce the 
following formula due to Chiswell [5] (our formula differs slightly from his, 
because he includes the empty subgraph and we do not.

\proclaim{Corollary 2} With the above notation, $\chi(G \Gamma) = \sum 
\chi_K$, the sum being taken over all complete subgraphs of $\Gamma$. 
\endproclaim

Let $\langle Y;S \rangle$ be a presentation of a group. This group then 
has a free resoltuion whose basis in dimension 1 is (a set bijective with) 
$Y$, and whose basis in dimension 2 is (a set bijective with) $S$. The 
boundary operator in dimension 2 can be described by means of the Fox 
derivatives.

Take such a presentation $\langle Y_v;S_v \rangle$ for each vertex group 
of a graph product, and take the corresponding free resolutions. Then the 
graph product has a presentation $\langle \bigcup Y_v; \bigcup S_v \cup C$, 
where $C$ is the set of all $[a,b]$ for $a \in Y_v, b \in Y_w$ and all 
adjacent vertices $v$ and $w$ such that $v<w$ in the chosen ordering of the 
graph.

Consider the resolution of the graph product obtained by our theorem from 
the free resolutions of the vertex groups. It will be a free resolution whose 
basis in dimension 1 will be $\bigcup Y_v$. In dimension 2 it will have two 
kinds of basis elements. The members of $\bigcup S_v$ form the first kind, 
while the second kind consists of all $a \otimes b$, where $a \in Y_v, b \in 
Y_w$, for all adjacent vertices $v$ and $w$ such that $v<w$. Thus this basis 
is bijective with the basis given by the presentation of the graph product. 
Further, the boundary operator given by the theorem is exactly that given by 
the presentation.

In dimension 3 there are four kinds of basis elements. The first kind 
consists of all basis elements in dimension 3 of the resolutions of the 
vertex groups. The second kind consists of all $a \otimes s$ where $a \in 
Y_v, s \in S_w$, for all adjacent vertices $v$ and $w$ such that $v<w$. The 
third kind consists of all $s \otimes b$ where $s \in S_v, b\in Y_w$, for all 
adjacent vertices $v$ and $w$ such that $v<w$. The fourth kind consists of 
all $a \otimes b \otimes c$ where $a \in Y_u, b \in Y_v, c \in Y_w$ and $\{ 
u,v,w \}$ is a complete subgraph with $u<v<w$.

The boundary operator can be explicitly described for the basis elements 
of the second, third, and fourth kinds. When this is done, we recover a 
result of Baik, Howie, and Pride [1].


\Refs
\ref \no 1 \by Y-G. Baik, J.Howie and S.J. Pride \paper The
identity problem for graph products \jour J. Algebra
\toappear \endref
\ref \no 2 \by H. Bass \paper Euler characteristics and
characters of discrete groups \jour Invent. Math. \vol 35
\yr 1976 \pages 155--196 \endref
\ref no 3 \by K.S. Brown \book Cohomology of groups \bookinfo Graduate
texts in mathematics {\bf 87} \publ Springer-Verlag \yr 1982 \endref
\ref \no 4 \by I.M.Chiswell \paper Euler characteristics
of groups \jour Math. Zeit. \vol 147 \yr 1976 \pages 1-11
\endref
\ref \no 5 \by I.M. Chiswell \paper The Euler
characteristic of a graph product \endref
\ref \no 6 \by D.E. Cohen \book Groups of cohomolgical
dimension one \bookinfo Springer lecture notes in mathematics \publ
Springer-Verlag \yr 1972 \endref \ref \no 7 \by E.R. Green \book Graph products
of groups \bookinfo Ph.D. Thesis \publ University of Leeds \yr 1990 \endref \ref
\no 8 \by A. Hattori \paper Rank element of a projective module \jour Nagoya
Math. J. \vol 25 \yr 1965 \pages 113--120 \endref
\ref \no 9 \by J.R. Stallings \paper Centerless groups
--- an algebraic formulation of Gottlieb's theorem \jour
Topology \vol 4 \yr 1965 \pages 129--134 \endref
\endRefs
 \enddocument